\long\def\symbolfootnote[#1]#2{\begingroup%
\def\thefootnote{\fnsymbol{footnote}}\footnote[#1]{#2}\endgroup}

\documentclass[12pt]{article}
\usepackage{amsfonts}

\usepackage{amsmath}

\newtheorem{theorem}{Theorem}

\newtheorem{lemma}{Lemma}

\begin{document}

\title{A sharp estimate for the bottom of the spectrum of the Laplacian on
K\"{a}hler manifolds}
\author{Ovidiu Munteanu}
\date{}
\maketitle

\begin{abstract}
On a complete noncompact K\"{a}hler manifold we prove that the bottom of the
spectrum for the Laplacian is bounded from above by $m^2$ if the Ricci
curvature is bounded from below by $-2(m+1)$. Then we show that if this
upper bound is achieved then the manifold has at most two ends. These results
improve previous results on this subject proved by P. Li and J. Wang in \cite
{L-W3} and \cite{L-W} under assumptions on the bisectional curvature.
\end{abstract}

\overfullrule=5pt

\section{Introduction}

\symbolfootnote[0]{
Reaserch partially supported by NSF grant No. DMS-0503735}

Let $N^n$ be a complete noncompact Riemannian manifold of dimension $n$ and
assume that the Ricci curvature has the lower bound $Ric\geq -(n-1).$ As
a consequence of Cheng's theorem (\cite{C}) we know that if $\lambda
_1\left( N\right) $ denotes the bottom of the spectrum of the Laplacian on $N
$ then $\lambda _1\left( N\right) \leq \frac{\left( n-1\right) ^2}4.$ This
is a sharp upper bound for $\lambda _1\left( N\right) $, the hyperbolic
space form $\mathbb{H}^n$ is an example where equality is achieved. Recall that
the proof of Cheng's theorem is relying on the Laplacian comparison theorem
for Riemannian manifolds, that is to say on an upper bound of the Laplacian
of the distance function on $N$. An interesting question is to study all
manifolds satisfying the equality case in Cheng's theorem, i.e. those
manifolds for which $Ric\geq -(n-1)$ and $\lambda _1\left( N\right) =\frac{%
\left( n-1\right) ^2}4.$ In \cite{L-W2} P. Li and J. Wang proved that if
equality holds in Cheng's upper bound and $n\geq 3$ then either the manifold
has one end or the manifold has two ends in which case $N$ must either be

(1) a warped product $N=\mathbb{R}\times P$ with $P$ compact and metric given
by $ds_N^2=dt^2+\exp \left( 2t\right) ds_P^2$, or

(2) if $n=3$ a warped product $N=\mathbb{R}\times P$ with $P$ compact and
metric given by $ds_N^2=dt^2+\cosh ^2\left( t\right) ds_P^2$.

In \cite{L-W3} P. Li and J. Wang have proved that Cheng's theorem has an
analogue in the K\"{a}hler setting. Throughout this paper $M^m$ is
a complete noncompact K\"{a}hler manifold of complex dimension $m.$

If $ds^2=h_{\alpha \overline{\beta }}dz^\alpha d\overline{z}^\beta $ denotes
the K\"{a}hler metric on $M,$ then ${Re}\left( ds^2\right) $ defines a
Riemannian metric on $M.$ Suppose $\left\{ e_1,e_2,...,e_{2m}\right\} $ with 
$e_{2k}=Je_{2k-1}$ for any $k\in \left\{ 1,...,m\right\} $ is an orthonormal
frame with respect to this Riemannian metric, then $\left\{
v_1,...,v_m\right\} $ is a unitary frame of $T_x^{1,0}M,$ where 
\[
v_k=\frac 12(e_{2k-1}-\sqrt{-1}e_{2k}).
\]
Recall that the bisectional curvature $BK_M$ of $M$ is defined by 
\[
R_{\alpha \bar{\alpha}\beta \bar{\beta}}=<R_{v_\alpha v_{\bar{\alpha}%
}}v_\beta ,v_{\bar{\beta}}>
\]
and we say that $BK_M\geq -1$ on $M$ if for any $\alpha$ and $\beta$
\[
R_{\alpha \bar{\alpha}\beta \bar{\beta}}\geq -(1+\delta _{\alpha \bar{\beta}%
}).
\]
Note that for the space form $\mathbb{C}\mathbb{H}^m$ we have $BK_{\mathbb{C}\mathbb{H}%
^m}=-1.$

\begin{theorem}
(P. Li and J. Wang) If $M^m$ is a complete noncompact K\"{a}hler manifold 
of complex dimension $m$ with $BK_M\geq -1$ then 
\[
\lambda _1\left( M\right) \leq m^2=\lambda _1\left( \mathbb{C}\mathbb{H}^m\right) .
\]
\end{theorem}

Li-Wang proved this theorem in the spirit of Cheng's proof, they first
obtained a Laplacian comparison theorem for manifolds with $BK_M\geq -1$
(Theorem 1.6 in \cite{L-W3}) and then the sharp estimate for $\lambda
_1\left( M\right) $ follows. We would like to point out that the bisectional
curvature assumption is essential in their proof of the Laplacian comparison
theorem.

An interesting question that one can ask is if the sharp estimate for 
$\lambda _1\left( M\right)$ from Theorem 1 remains true under a Ricci curvature
bound from below. 
This question is motivated in part by the following situation in the compact 
K\"{a}hler case, where we have a version of Lichnerowicz's theorem. Namely,
if for a compact K\"{a}hler manifold $N^m$ the Ricci curvature has
the lower bound $Ric_N\geq 2(m+1)$, then the first eigenvalue of the Laplacian
has a sharp lower bound, $\lambda_1\left( N\right)\geq 4(m+1)$. We are grateful to 
Lei Ni for pointing out this result to us, for a simple proof of it see \cite{U}.

In this paper, our
first goal is to show that indeed there is a sharp estimate for $\lambda _1\left( M\right)$
under only Ricci curvature lower bound. Our proof is based
on the variational principle for $\lambda _1\left( M\right) $ and
integration by parts. In fact, our argument can be localized on each end of
the manifold. 

\begin{theorem}
Let $M^m$ be a complete noncompact K\"{a}hler manifold of complex dimension $%
m$ such that the Ricci curvature is bounded from below by 
\[
Ric\geq -2\left( m+1\right) .
\]
Then if $E$ is an end of $M$ and $\lambda _1\left( E\right) $ is the infimum
of the Dirichlet spectrum of the Laplacian on $E\,$then 
\[
\lambda _1\left( E\right) \leq m^2.
\]
In particular, we have the sharp estimate 
\[
\lambda _1\left( M\right) \leq m^2.
\]
\end{theorem}
Note that the condition on the Ricci curvature in Theorem 2 means 
\[
Ric(e_k,e_j)\geq -2\left( m+1\right) \delta _{kj}
\]
for any $k,j\in \{1,..,2m\},$ which is equivalent to 
\[
Ric_{\alpha \overline{\beta }}\geq -\left( m+1\right) \delta _{\alpha 
\overline{\beta }},
\]
for the unitary frame $\left\{ v_1,v_2,...,v_m\right\} .$

The same as in the Riemannian setting, it is interesting to study the
K\"{a}hler manifolds for which equality is achieved in Theorem 2. Let us
recall that for bisectional curvature lower bound Li-Wang have proved in 
\cite{L-W} that such manifolds need to have at most two ends.

\begin{theorem}
(P. Li and J. Wang) If $M^m$ is a complete noncompact K\"{a}hler manifold with $%
\lambda _1\left( M\right) =m^2$ and $BK_M\geq -1$ then $M$ has at most two ends.
\end{theorem}

Their proof relies on a study of the Buseman function $\beta $ on $M,$ so the 
Laplacian comparison theorem plays again an important role. An intersting fact about 
their proof is that it gives an unified approach of the question in the Riemannian
and K\"{a}hlerian case.

We now want to make some comments on the case when $M$ has exactly two ends. 
In this case, the proof of Li-Wang provides some structure information of the manifold. 
Namely, not only that for any $t$ the level set $\beta=t$ is diffeomorphic to
the level set $\beta=t_0$ for some $t_0$ fixed, but also the metric on $\beta=t$ 
is determined by the metric on $\beta=t_0$. 

Our second goal in this paper is to obtain the same conclusion on the number of ends
if equality is achieved in Theorem 2. This will be done by a
careful study of the estimates in Theorem 2. 

If $M$ is assumed to have exactly two ends, using our approach we will be able to 
deduce the same structure information of the manifold as discussed above, for the level sets 
of a function defined by the Li-Tam theory.

\textbf{Remark}. 
After this paper was written the author was informed by Peter Li that the analysis of the two 
ends case can be deepened and in fact if $M$ has bounded curvature, then it is
isometrically covered by $\mathbb{C}\mathbb{H}^m$. Examples are also known, with both bounded and 
unbounded curvature, see \cite{L-W}. We expect this result can be recovered with our way, 
and in fact this will become apparent towards the end of the proof of Theorem 4.

\begin{theorem}
Let $M^m$ be a complete noncompact K\"{a}hler manifold of complex dimension $%
m$ such that the Ricci curvature is bounded from below by 
\[
Ric\geq -2\left( m+1\right) .
\]
If $\lambda _1\left( M\right) =m^2$ then $M$ has at most two ends.
\end{theorem}

\textbf{Aknowledgement}. The author would like to express his deep gratitude
to his advisor, Professor Peter Li, for constant help, support and many
valuable discussions.

\section{The proofs}

To prove Theorem 2 and Theorem 4 we first need the following preparation.

Let $E$ be a nonparabolic end of $M.$

Withouth loss of generality we will henceforth assume that $\lambda _1\left(
E\right) >0.$

From the theory of Li-Tam (\cite{L-T}) we know that there exists a harmonic
function $f$ on $E$ that is obtained with the following procedure. Let $f_R$
be the harmonic function with Dirichlet boundary conditions: $f_R=1$ on $%
\partial E,$ $f_R=0$ on $\partial E_p(R),$ where $E_p(R)=E\cap B_p\left(
R\right) .$

Then $f_R$ admits a subsequence convergent to $f,$ with the properties: $%
0<f<1$ on $E,$ $f$ $=1$ on $\partial E$ and $f$ has finite Dirichlet
integral. Moreover, since $\lambda _1\left( E\right) >0,$ we know by a
theorem of Li-Wang that (\cite{L-W1}) 
\[
\int_{E_p\left( R+1\right) \backslash E_p\left( R\right) }f^2\leq c_1\exp
\left( -2\sqrt{\lambda _1(E)}R\right) . 
\]
Further on integration on the level sets of $f$ will play a central role in
our proofs, and for this let us recall the following important property of $%
f $ (\cite{L-W4}). Namely, for $t,a,b<1$ let 
\[
l\left( t\right) =\left\{ x\in E\left| \;f\left( x\right) =t\right. \right\} 
\]
and define the set 
\[
L\left( a,b\right) =\left\{ x\in E\left| \;a<f\left( x\right) <b\right.
\right\} . 
\]
Then for almost all $t<1$ 
\[
\int_{l\left( t\right) }\left| \nabla f\right| =const<\infty 
\]
and we have: 
\[
\int_{L\left( a,b\right) }\left| \nabla f\right| ^2=\left( b-a\right)
\int_{l\left( t_0\right) }\left| \nabla f\right| . 
\]
Let us denote 
\[
L=L(\frac 12\delta \varepsilon ,2\varepsilon ), 
\]
where $\delta ,\varepsilon >0$ are sufficiently small fixed numbers to be
chosen later.

Since we will often integrate by parts on $L$, let us construct a cut-off $%
\phi $ with compact support in $L.$ Define $\phi =\psi \varphi $ with $\psi $
depending on the distance function

\[
\psi =\left\{ 
\begin{array}{c}
1 \\ 
R-r \\ 
0
\end{array}
\left. 
\begin{array}{c}
\text{on} \\ 
\text{on} \\ 
\text{on}
\end{array}
\right. \left. 
\begin{array}{l}
E_p\left( R-1\right) \\ 
E_p\left( R\right) \backslash E_p\left( R-1\right) \\ 
E\backslash E_p\left( R\right)
\end{array}
\right. \right. 
\]
and $\varphi $ defined on the level sets of $f$%
\[
\varphi =\left\{ 
\begin{array}{c}
(\log 2)^{-1}(\log f-\log (\frac 12\delta \varepsilon )) \\ 
1 \\ 
(\log 2)^{-1}(\log 2\varepsilon -\log f) \\ 
0
\end{array}
\left. 
\begin{array}{l}
\text{on}\;\;L(\frac 12\delta \varepsilon ,\delta \varepsilon ) \\ 
\text{on}\;L\left( \delta \varepsilon ,\varepsilon \right) \\ 
\text{on}\;\;L\left( \varepsilon ,2\varepsilon \right) \\ 
\text{otherwise.}
\end{array}
\right. \right. 
\]
For convenience, let us assume $R=\frac 1{\delta \varepsilon }.$ We have the
following result:

\begin{lemma}
For any $0<a<2$ the following inequality holds: 
\begin{eqnarray*}
\frac 1{16}\left( \frac 1m-\frac{\left( 1-a\right) ^2}{a\left( 2-a\right) }%
\right) \frac 1{(-\log \delta )}\int_L\frac{\left| \nabla f\right| ^4}{f^3}%
\phi ^2 &\leq &\frac a{2-a}\frac{m+1}4\int_{l\left( t_0\right) }\left|
\nabla f\right|  \\
&&+\frac c{(-\log \delta )^{\frac 12}},
\end{eqnarray*}
where $c$ is a constant not depending on $\delta $ or $\varepsilon .$
\end{lemma}
\textbf{Proof of Lemma 1.} Note that the gradient and the Laplacian satisfy: 
\begin{eqnarray*}
\nabla f\cdot \nabla g &=&2\left( f_\alpha g_{\overline{\alpha }}+f_{%
\overline{\alpha }}g_\alpha \right) \\
\Delta f &=&4f_{\alpha \overline{\alpha }}.
\end{eqnarray*}

Let $u=\log f,$ then a simple computation shows that 
\[
u_{\alpha \overline{\beta }}=f^{-1}f_{\alpha \overline{\beta }%
}-f^{-2}f_\alpha f_{\overline{\beta }}. 
\]

Consider now 
\[
\int_Lf\left| u_{\alpha \overline{\beta }}\right| ^2\phi ^2 
\]
which we estimate from above and from below to prove our claim.

To begin with, 
\[
\int_Lf\left| u_{\alpha \overline{\beta }}\right| ^2\phi
^2=\int_Lf^{-1}\left| f_{\alpha \overline{\beta }}\right| ^2\phi
^2-2\int_Lf^{-2}(f_{\alpha \overline{\beta }}f_{\overline{\alpha }}f_\beta
)\phi ^2+\frac 1{16}\int_Lf^{-3}\left| \nabla f\right| ^4\phi ^2. 
\]
The first term is 
\begin{eqnarray*}
\int_Lf^{-1}\left| f_{\alpha \overline{\beta }}\right| ^2\phi ^2
&=&\int_Lf^{-1}(f_{\alpha \overline{\beta }}\cdot f_{\overline{\alpha }\beta
})\phi ^2=-\int_Lf_\alpha \left( f^{-1}f_{\overline{\alpha }\beta }\phi
^2\right) _{\overline{\beta }} \\
&=&\int_Lf^{-2}(f_{\overline{\alpha }\beta }f_\alpha f_{\overline{\beta }%
})\phi ^2-\int_Lf^{-1}f_\alpha f_{\overline{\alpha }\beta \overline{\beta }%
}\phi ^2-\int_Lf^{-1}f_{\overline{\alpha }\beta }f_\alpha \left( \phi
^2\right) _{\overline{\beta }}
\end{eqnarray*}
and using the Ricci identities and $\Delta f=0$ we see that $f_{\overline{%
\alpha }\beta \overline{\beta }}=0.$ It also shows that the last integral
needs to be a real number.

This proves that

\begin{eqnarray}
\int_Lf\left| u_{\alpha \overline{\beta }}\right| ^2\phi ^2
&=&-\int_Lf^{-2}(f_{\alpha \overline{\beta }}f_{\overline{\alpha }}f_\beta
)\phi ^2  \nonumber \\
&&+\frac 1{16}\int_Lf^{-3}\left| \nabla f\right| ^4\phi ^2-\int_Lf^{-1}f_{%
\overline{\alpha }\beta }f_\alpha \left( \phi ^2\right) _{\overline{\beta }}.
\label{1}
\end{eqnarray}
Let us use again integration by parts to see that 
\begin{eqnarray*}
-\int_Lf^{-2}(f_{\alpha \overline{\beta }}f_{\overline{\alpha }}f_\beta
)\phi ^2 &=&\int_Lf_\alpha \left( f^{-2}f_{\overline{\alpha }}f_\beta \phi
^2\right) _{\overline{\beta }}= \\
&=&-2\int_Lf^{-3}f_\alpha f_{\overline{\alpha }}f_\beta f_{\overline{\beta }%
}\phi ^2+\int_Lf^{-2}f_{\overline{\alpha }\overline{\beta }}f_\alpha f_\beta
\phi ^2 \\
&&+\int_Lf^{-2}f_\alpha f_{\overline{\alpha }}f_\beta \left( \phi ^2\right)
_{\overline{\beta }}.
\end{eqnarray*}
Similarly, one finds 
\begin{eqnarray*}
-\int_Lf^{-2}(f_{\alpha \overline{\beta }}f_{\overline{\alpha }}f_\beta
)\phi ^2 &=&-\int_Lf^{-2}(f_{\overline{\alpha }\beta }f_\alpha f_{\bar{\beta}%
})\phi ^2=\int_Lf_{\bar{\alpha}}\left( f^{-2}f_\alpha f_{\bar{\beta}}\phi
^2\right) _\beta \\
&=&-2\int_Lf^{-3}f_\alpha f_{\overline{\alpha }}f_\beta f_{\overline{\beta }%
}\phi ^2+\int_Lf^{-2}f_{\alpha \beta }f_{\bar{\alpha}}f_{\bar{\beta}}\phi ^2
\\
&&+\int_Lf^{-2}f_\alpha f_{\overline{\alpha }}f_{\bar{\beta}}\left( \phi
^2\right) _\beta .
\end{eqnarray*}
Combining the two identities we get 
\begin{eqnarray}
-\int_Lf^{-2}(f_{\alpha \overline{\beta }}f_{\overline{\alpha }}f_\beta
)\phi ^2 &=&-\frac 18\int_Lf^{-3}\left| \nabla f\right| ^4\phi
^2+\int_Lf^{-2}{Re}\left( f_{\overline{\alpha }\overline{\beta }}f_\alpha
f_\beta \right)  \nonumber \\
&&+\frac 14\int_Lf^{-2}\left| \nabla f\right| ^2Re(f_{\bar{\beta}}\left(
\phi ^2\right) _\beta ).  \label{2}
\end{eqnarray}
Note that the following inequality holds on $E$: 
\begin{equation}
\left| f_{\overline{\alpha }\overline{\beta }}f_\alpha f_\beta \right| \leq
\frac 14\left| f_{\alpha \beta }\right| \left| \nabla f\right| ^2  \label{3}
\end{equation}

We want to insist on the proof of this inequality because it will matter
when we study the manifolds with $\lambda _1\left( M\right) =m^2$. Since the
two numbers in (\ref{3}) are independent of the unitary frame, let us choose
a frame at the fixed point $x\in E$ such that 
\[
e_1=\frac 1{\left| \nabla f\right| }\nabla f. 
\]
Certainly, we need $\left| \nabla f\right| \left( x\right) \neq 0$ which we
assume without loss of generality because if $\left| \nabla f\right| \left(
x\right) =0$ there is nothing to prove.

Then one can see that 
\[
f_{e_1}=\left| \nabla f\right| ,\;f_{e_2}=0,....f_{e_{2m}}=0 
\]
or, in the unitary frame 
\[
f_1=f_{\bar{1}}=\frac 12\left| \nabla f\right| ,\;\;\;f_\alpha =f_{\bar{%
\alpha}}=0\;\;\text{if}\;\alpha >1. 
\]
This proves the inequality because 
\[
\left| f_{\overline{\alpha }\overline{\beta }}f_\alpha f_\beta \right|
=\frac 14\left| \nabla f\right| ^2\left| f_{11}\right| \leq \frac 14\left|
f_{\alpha \beta }\right| \left| \nabla f\right| ^2. 
\]
Moreover, we learn that equality holds in (\ref{3}) if and only if 
\[
f_{\alpha \beta }=0\;\;\text{for}\;(\alpha ,\beta )\neq \left( 1,1\right) , 
\]
with respect to the frame chosen above.

Since the following holds: 
\[
Re\left( f_{\overline{\alpha }\overline{\beta }}f_\alpha f_\beta \right)
\leq \left| f_{\overline{\alpha }\overline{\beta }}f_\alpha f_\beta \right|
\leq \frac 14\left| f_{\alpha \beta }\right| \left| \nabla f\right| ^2, 
\]
we get for an arbitrary $a>0$ 
\begin{gather}
2\int_Lf^{-2}Re\left( f_{\overline{\alpha }\overline{\beta }}f_\alpha
f_\beta \right) \phi ^2\leq \int_L2\left( f^{-1/2}\left| f_{\alpha \beta
}\right| \phi \right) \left( \frac 14f^{-3/2}\left| \nabla f\right| ^2\phi
\right)  \nonumber \\
\leq a\int_Lf^{-1}\left| f_{\alpha \beta }\right| ^2\phi ^2+\frac
1{16a}\int_Lf^{-3}\left| \nabla f\right| ^4\phi ^2.  \label{4}
\end{gather}
Moreover, again integrating by parts we have 
\begin{eqnarray*}
\int_Lf^{-1}\left| f_{\alpha \beta }\right| ^2\phi ^2
&=&\int_Lf^{-1}f_{\alpha \beta }f_{\bar{\alpha}\bar{\beta}}\phi
^2=-\int_Lf_\alpha \left( f^{-1}f_{\bar{\alpha}\bar{\beta}}\phi ^2\right)
_\beta \\
&=&\int_Lf^{-2}f_{\overline{\alpha }\overline{\beta }}f_\alpha f_\beta \phi
^2-\int_Lf^{-1}f_\alpha f_{\bar{\alpha}\bar{\beta}\beta }\phi
^2-\int_Lf^{-1}f_\alpha f_{\bar{\alpha}\bar{\beta}}\left( \phi ^2\right)
_\beta
\end{eqnarray*}
and on the other hand 
\begin{eqnarray*}
\int_Lf^{-1}\left| f_{\alpha \beta }\right| ^2\phi ^2
&=&\int_Lf^{-1}f_{\alpha \beta }f_{\bar{\alpha}\bar{\beta}}\phi ^2=-\int_Lf_{%
\bar{\alpha}}\left( f^{-1}f_{\alpha \beta }\phi ^2\right) _{\bar{\beta}} \\
&=&\int_Lf^{-2}f_{\alpha \beta }f_{\bar{\alpha}}f_{\bar{\beta}}\phi
^2-\int_Lf^{-1}f_{\bar{\alpha}}f_{\alpha \beta \bar{\beta}}\phi
^2-\int_Lf^{-1}f_{\bar{\alpha}}f_{\alpha \beta }\left( \phi ^2\right) _{\bar{%
\beta}}
\end{eqnarray*}
so that combining the two identities we get 
\begin{eqnarray*}
\int_Lf^{-1}\left| f_{\alpha \beta }\right| ^2\phi ^2 &=&\int_Lf^{-2}Re(f_{%
\overline{\alpha }\overline{\beta }}f_\alpha f_\beta )\phi
^2-\int_Lf^{-1}f_\alpha f_{\bar{\alpha}\bar{\beta}\beta }\phi ^2 \\
&&-\int_Lf^{-1}Re(f_\alpha f_{\bar{\alpha}\bar{\beta}}\left( \phi ^2\right)
_\beta ).
\end{eqnarray*}
Note that the Ricci identities imply 
\[
f_{\bar{\alpha}\bar{\beta}\beta }=f_{\bar{\beta}\bar{\alpha}\beta }=f_{\bar{%
\beta}\bar{\beta}\bar{\alpha}}+Ric_{\beta \bar{\alpha}}f_{\bar{\beta}} 
\]
and therefore we obtain 
\begin{eqnarray*}
\int_Lf^{-1}\left| f_{\alpha \beta }\right| ^2\phi ^2 &\leq
&\int_Lf^{-2}Re(f_{\overline{\alpha }\overline{\beta }}f_\alpha f_\beta
)\phi ^2+\frac{m+1}4\int_Lf^{-1}\left| \nabla f\right| ^2\phi ^2 \\
&&-\int_Lf^{-1}Re(f_\alpha f_{\bar{\alpha}\bar{\beta}}\left( \phi ^2\right)
_\beta ).
\end{eqnarray*}
Plug this inequality into (\ref{4}) and it follows 
\begin{gather*}
\left( 2-a\right) \int_Lf^{-2}Re\left( f_{\overline{\alpha }\overline{\beta }%
}f_\alpha f_\beta \right) \phi ^2\leq a\frac{m+1}4\int_Lf^{-1}\left| \nabla
f\right| ^2\phi ^2 \\
+\frac 1{16a}\int_Lf^{-3}\left| \nabla f\right| ^4\phi
^2-a\int_Lf^{-1}Re(f_\alpha f_{\bar{\alpha}\bar{\beta}}\left( \phi ^2\right)
_\beta ).
\end{gather*}
Let us fix henceforth $0<a<2$ to make sure that $2-a>0.$ Now we are getting
back to (\ref{2}) and obtain 
\begin{gather}
-\int_Lf^{-2}(f_{\alpha \overline{\beta }}f_{\overline{\alpha }}f_\beta
)\phi ^2\leq \left( -\frac 18+\frac 1{16a(2-a)}\right) \int_Lf^{-3}\left|
\nabla f\right| ^4\phi ^2  \nonumber \\
+\frac a{2-a}\frac{m+1}4\int_Lf^{-1}\left| \nabla f\right| ^2\phi ^2-\frac
a{2-a}\int_Lf^{-1}Re(f_\alpha f_{\bar{\alpha}\bar{\beta}}\left( \phi
^2\right) _\beta )  \nonumber \\
+\frac 14\int_Lf^{-2}\left| \nabla f\right| ^2Re(f_{\bar{\beta}}\left( \phi
^2\right) _\beta ).  \label{5}
\end{gather}
We have thus proved that 
\begin{gather}
\int_Lf\left| u_{\alpha \overline{\beta }}\right| ^2\phi ^2\leq \left(
-\frac 1{16}+\frac 1{16a(2-a)}\right) \int_Lf^{-3}\left| \nabla f\right|
^4\phi ^2  \nonumber \\
+\frac a{2-a}\frac{m+1}4\int_Lf^{-1}\left| \nabla f\right| ^2\phi ^2+\frac
14\int_Lf^{-2}\left| \nabla f\right| ^2Re(f_{\bar{\beta}}\left( \phi
^2\right) _\beta )  \nonumber \\
-\int_Lf^{-1}f_{\overline{\alpha }\beta }f_\alpha \left( \phi ^2\right) _{%
\overline{\beta }}-\frac a{2-a}\int_Lf^{-1}Re(f_\alpha f_{\bar{\alpha}\bar{%
\beta}}\left( \phi ^2\right) _\beta ).  \label{6}
\end{gather}
To finish the upper estimate of $\int_Lf\left| u_{\alpha \overline{\beta }%
}\right| ^2\phi ^2$ we need to estimate the terms involving $\left( \phi
^2\right) _\beta .$ We will prove that they can be bounded from above by a
constant $\cdot (-\log \delta )^{1/2}.$

Start with 
\begin{eqnarray*}
2\int_Lf^{-2}\left| \nabla f\right| ^2Re(f_{\bar{\beta}}\left( \phi
^2\right) _\beta ) &\leq &\frac 12\int_Lf^{-2}\left| \nabla f\right|
^3\left| \nabla \phi ^2\right| \\
&\leq &\int_Lf^{-2}\left| \nabla f\right| ^3\left| \nabla \varphi \right|
\psi +\int_Lf^{-2}\left| \nabla f\right| ^3\left| \nabla \psi \right| \varphi
\end{eqnarray*}
Now it is easy to see that by the gradient estimate and co-area formula 
\begin{eqnarray*}
\int_Lf^{-2}\left| \nabla f\right| ^3\left| \nabla \varphi \right| &\leq
&c_2(\int_{L(\frac 12\delta \varepsilon ,\delta \varepsilon )}f^{-1}\left|
\nabla f\right| ^2+\int_{L(\varepsilon ,2\varepsilon )}f^{-1}\left| \nabla
f\right| ^2) \\
&\leq &c_{3,}
\end{eqnarray*}
while by the decay rate of $f^2$ we get 
\[
\int_Lf^{-2}\left| \nabla f\right| ^3\left| \nabla \psi \right| \leq
c_4\frac 1{\delta \varepsilon }\exp \left( -2\sqrt{\lambda _1(E)}R\right)
\leq c_5, 
\]
using that $R=\frac 1{\delta \varepsilon }.$ Clearly, the constants so far
do not depend on the choice of $\delta $ or $\varepsilon .$

To estimate the other terms one proceeds similarly. For example, 
\begin{eqnarray*}
-2\int_Lf^{-1}Re(f_\alpha f_{\bar{\alpha}\bar{\beta}}\left( \phi ^2\right)
_\beta ) &\leq &\int_Lf^{-1}\left| f_{\bar{\alpha}\bar{\beta}}\right| \left|
\nabla f\right| \phi \left| \nabla \phi \right| \\
&\leq &\left( \int_Lf^{-1}\left| \nabla f\right| ^2\left| \nabla \phi
\right| ^2\right) ^{\frac 12}\left( \int_Lf^{-1}\left| f_{\bar{\alpha}\bar{%
\beta}}\right| ^2\phi ^2\right) ^{\frac 12} \\
&\leq &c_6\left( \int_Lf^{-1}\left| f_{\bar{\alpha}\bar{\beta}}\right|
^2\phi ^2\right) ^{\frac 12},
\end{eqnarray*}
However, using an inequality proved above we get 
\begin{eqnarray*}
\int_Lf^{-1}\left| f_{\alpha \beta }\right| ^2\phi ^2 &\leq
&\int_Lf^{-2}Re(f_{\overline{\alpha }\overline{\beta }}f_\alpha f_\beta
)\phi ^2+\frac{m+1}4\int_Lf^{-1}\left| \nabla f\right| ^2\phi ^2 \\
&&-\int_Lf^{-1}Re(f_\alpha f_{\bar{\alpha}\bar{\beta}}\left( \phi ^2\right)
_\beta ) \\
&\leq &\frac 14\int_Lf^{-2}\left| f_{\alpha \beta }\right| \left| \nabla
f\right| ^2\phi ^2+\frac{m+1}4\int_Lf^{-1}\left| \nabla f\right| ^2\phi ^2 \\
&&+\frac 12\int_Lf^{-1}\left| f_{\alpha \beta }\right| \left| \nabla
f\right| \phi \left| \nabla \phi \right| \\
&\leq &\frac 18\int_Lf^{-1}\left| f_{\alpha \beta }\right| ^2\phi ^2+\frac
18\int_Lf^{-3}\left| \nabla f\right| ^4\phi ^2 \\
&&+\frac{m+1}4\int_Lf^{-1}\left| \nabla f\right| ^2\phi ^2 \\
&&+\frac 14\int_Lf^{-1}\left| f_{\alpha \beta }\right| ^2\phi ^2+\frac
14\int_Lf^{-1}\left| \nabla f\right| ^2\left| \nabla \phi \right| ^2,
\end{eqnarray*}
which shows there exists constants $c_7$ and $c_8$ such that: 
\begin{eqnarray*}
\int_Lf^{-1}\left| f_{\bar{\alpha}\bar{\beta}}\right| ^2\phi ^2 &\leq
&c_7\int_Lf^{-1}\left| \nabla f\right| ^2\phi ^2+c_8\int_Lf^{-1}\left|
\nabla f\right| ^2\left| \nabla \phi \right| ^2 \\
&\leq &c_9\left( -\log \delta \right) .
\end{eqnarray*}
We have proved that 
\[
\int_Lf^{-1}\left| f_{\bar{\alpha}\bar{\beta}}\right| \left| \nabla f\right|
\phi \left| \nabla \phi \right| \leq c_{10}(-\log \delta )^{\frac 12}. 
\]
Let us gather the information we have so far: 
\begin{eqnarray*}
\int_Lf\left| u_{\alpha \overline{\beta }}\right| ^2\phi ^2 &\leq &\left(
-\frac 1{16}+\frac 1{16a(2-a)}\right) \int_Lf^{-3}\left| \nabla f\right|
^4\phi ^2 \\
&&+\frac a{2-a}\frac{m+1}4\int_Lf^{-1}\left| \nabla f\right| ^2\phi
^2+c(-\log \delta )^{\frac 12}.
\end{eqnarray*}
The estimate from below is straightforward: 
\begin{equation}
\left| u_{\alpha \overline{\beta }}\right| ^2\geq \sum_\alpha \left|
u_{\alpha \bar{\alpha}}\right| ^2\geq \frac 1m\left| \sum_\alpha u_{\alpha 
\overline{\alpha }}\right| ^2=\frac 1{16m}f^{-4}\left| \nabla f\right| ^4.
\label{7}
\end{equation}
Hence, this shows that 
\begin{eqnarray*}
\frac 1{16}\left( 1+\frac 1m-\frac 1{a(2-a)}\right) \int_Lf^{-3}\left|
\nabla f\right| ^4\phi ^2 &\leq &\frac a{2-a}\frac{m+1}4\int_Lf^{-1}\left|
\nabla f\right| ^2\phi ^2 \\
&&+c(-\log \delta )^{\frac 12},
\end{eqnarray*}
which proves the Lemma. \textbf{Q.E.D.}

In the following Lemma, we will estimate $\int_Lf^{-3}\left| \nabla f\right|
^4\phi ^2$ from bellow. To serve our purpose, we need this estimate to
depend on $\lambda _1\left( E\right) $ and this is done using the
variational principle. Recall that $E$ is a nonparabolic end, $\lambda
_1\left( E\right) >0$ and we set $L=L\left( \frac 12\delta \varepsilon
,2\varepsilon \right) $ for $\delta ,\varepsilon $ sufficiently small.

\begin{lemma}
\[
\frac 1{(-\log \delta )}\int_Lf^{-3}\left| \nabla f\right| ^4\phi ^2\geq
4\lambda _1\left( E\right) \int_{l(t_0)}\left| \nabla f\right| -\frac{c_0}{%
(-\log \delta )^{\frac 12}}.
\]
\end{lemma}
\textbf{Proof of Lemma 2.}

By the variational principle for $\lambda _1\left( E\right) ,$ 
\[
\lambda _1\left( E\right) \int_Ef\phi ^2\leq \int_E\left| \nabla \left( \phi
f^{\frac 12}\right) \right| ^2, 
\]
which means that 
\begin{eqnarray*}
\lambda _1\left( E\right) \int_Lf\phi ^2 &\leq &\frac 14\int_Lf^{-1}\left|
\nabla f\right| ^2\phi ^2+\int_Lf\left| \nabla \phi \right| ^2+\int_L\phi
\left| \nabla f\right| \left| \nabla \phi \right| \\
&\leq &\frac 14\int_{L\left( \delta \varepsilon ,\varepsilon \right)
}f^{-1}\left| \nabla f\right| ^2+c_{11},
\end{eqnarray*}
based on estimates similar to what we did in Lemma 1.

This implies that 
\[
\frac 1{(-\log \delta )}\int_Lf\phi ^2\leq \frac 1{4\lambda _1\left(
E\right) }\int_{l(t_0)}\left| \nabla f\right| +\frac{c_{11}}{(-\log \delta )}%
. 
\]
Finally, using the Schwarz inequality we get 
\begin{eqnarray*}
\int_{l(t_0)}\left| \nabla f\right| &=&\frac 1{(-\log \delta )}\int_{L\left(
\delta \varepsilon ,\varepsilon \right) }f^{-1}\left| \nabla f\right| ^2 \\
&\leq &\frac 1{(-\log \delta )}\int_Lf^{-1}\left| \nabla f\right| ^2\phi
^2+\frac 1{(-\log \delta )}\int_{L\cap \left( E\backslash E_p\left(
R-1\right) \right) }f^{-1}\left| \nabla f\right| ^2 \\
&\leq &\left( \frac 1{(-\log \delta )}\int_Lf^{-3}\left| \nabla f\right|
^4\phi ^2\right) ^{\frac 12}\left( \frac 1{(-\log \delta )}\int_Lf\phi
^2\right) ^{\frac 12} \\
&&+\frac{c_{12}}{(-\log \delta )}\frac 1{\delta \varepsilon }\exp \left( -2%
\sqrt{\lambda _1(E)}R\right) \\
&\leq &\left( \frac 1{(-\log \delta )}\int_Lf^{-3}\left| \nabla f\right|
^4\phi ^2\right) ^{\frac 12}\times \\
&&\times \left( \frac 1{4\lambda _1\left( E\right) }\int_{l(t_0)}\left|
\nabla f\right| +\frac{c_{11}}{(-\log \delta )}\right) ^{\frac 12}+\frac{%
c_{13}}{(-\log \delta )},
\end{eqnarray*}
which proves the Lemma. \textbf{Q.E.D}.

Suppose now that $M$ has a parabolic end $F.$ A theorem of Nakai (\cite{N},
see also \cite{N-R}) states that there exists an exhaustion function $f$ on $%
\overline{F}$ which is harmonic on $F$ and $f=0$ on $\partial F.$ In this
case we consider for $T,\beta >0$ fixed
\[
\phi =\left\{ 
\begin{array}{c}
(\log 2)^{-1}(\log f-\log (\frac 12T)) \\ 
1 \\ 
(\log 2)^{-1}(\log (2\beta T)-\log f) \\ 
0
\end{array}
\left. 
\begin{array}{l}
\text{on}\;\;L(\frac 12T,T) \\ 
\text{on}\;L\left( T,\beta T\right)  \\ 
\text{on}\;\;L\left( \beta T,2\beta T\right)  \\ 
\text{otherwise,}
\end{array}
\right. \right. 
\]
where the level sets are now defined on $F.$ Since $f$ is proper, there is
no need for a cut-off depending on the distance function. Our point now is
that Lemma 1 and Lemma 2 hold for this choice of $\phi $ also, the proofs
are identical. Note that if 
\[
\tilde{L}=L(\frac 12T,2\beta T)
\]
then the following inequalities hold on $\tilde{L}$: 
\begin{eqnarray*}
\frac 1{16}\left( \frac 1m-\frac{\left( 1-a\right) ^2}{a\left( 2-a\right) }%
\right) \frac 1{\log \beta }\int_{\tilde{L}}\frac{\left| \nabla f\right| ^4}{%
f^3}\phi ^2 &\leq &\frac a{2-a}\frac{m+1}4\int_{l\left( t_0\right) }\left|
\nabla f\right|  \\
&&+\frac{\widetilde{c}}{(\log \beta )^{\frac 12}},
\end{eqnarray*}
and 
\[
\frac 1{\log \beta }\int_{\tilde{L}}f^{-3}\left| \nabla f\right| ^4\phi
^2\geq 4\lambda _1\left( F\right) \int_{l(t_0)}\left| \nabla f\right| -\frac{%
\widetilde{c}_0}{(\log \beta )^{\frac 12}}.
\]

\bigskip

Now we are ready to prove Theorem 2.

\vspace{0.5in} \textbf{Proof of Theorem 2.}

Let us first prove the Theorem for a nonparabolic end $E.$

We know from Lemma 1 and Lemma 2 that 
\[
\lambda _1\left( E\right) \frac 14\left( \frac 1m-\frac{\left( 1-a\right) ^2%
}{a\left( 2-a\right) }\right) \int_{l(t_0)}\left| \nabla f\right| \leq \frac
a{2-a}\frac{m+1}4\int_{l\left( t_0\right) }\left| \nabla f\right| +\frac
C{(-\log \delta )^{\frac 12}}, 
\]
inequality that holds for any $\delta >0$ and for any $0<a<2.$

Therefore, making $\delta \rightarrow 0$ we get that for any $0<a<2,$%
\[
\lambda _1\left( E\right) \leq \frac{a\left( m+1\right) }{2-a}\left( \frac
1m-\frac{\left( 1-a\right) ^2}{a(2-a)}\right) ^{-1}. 
\]
Let us choose $a=\frac m{m+1},$ then 
\[
\frac 1m-\frac{\left( 1-a\right) ^2}{a(2-a)}=\frac 1m\frac{m+1}{m+2}\;;\;%
\frac{a\left( m+1\right) }{2-a}=\frac{m(m+1)}{m+2}, 
\]
which shows that

\[
\lambda _1\left( E\right) \leq m^2.
\]
This proves Theorem 2 for a nonparabolic end $E$. The proof for a parabolic
end $F$ is verbatim. \textbf{Q.E.D.}

\vspace{0.5in}

\textbf{Proof of Theorem 4:}

Suppose that $M$ has more than one end.
We know from \cite{L-W}
that if $\lambda _1\left( M\right) >\frac{m+1}2$ the manifold has only one
nonparabolic end. Hence let us set $E$ this nonparabolic end and
consequently $F=M\backslash E$ will be a parabolic end. Note that an end of $M$
is defined with respect to a compact subset of $M$, so that writing $M=E\cup F$ 
with $E$ nonparabolic and $F$ parabolic we are not loosing generality, in fact
$M$ can have many ends with respect to other compact subsets.   
The construction of
Li-Tam implies that there exists a harmonic function $f:M\rightarrow \left(
0,\infty \right) $ with the following properties:

1. On $E$ the function has the decay rate 
\[
\int_{E_p\left( R\right) \backslash E_p\left( R-1\right) }f^2\leq c_1\exp
\left( -2\sqrt{\lambda _1(M)}R\right) , 
\]

2. On $F$ the function is proper.

3. We have: 
\[
\sup_{x\in F}f\left( x\right) =\infty ,\;\;\inf_{x\in E}f\left( x\right) =0. 
\]

Let us point out some facts about the proofs of Lemma 1 and Lemma 2. In the
two lemmata, the function $f$ was defined only on a single end, which was
first assumed to be nonparabolic, and then we observed that  the proofs
still work on a parabolic end. In the framework of Theorem 4, we know that $f
$ is defined on the whole manifold, so now $L=L\left( b_0,b_1\right)
=\left\{ x\in M\left| \;b_0<f\left( x\right) <b_1\right. \right\} $ makes
sense for any $0<b_0<b_1$. One can see that the computations proved in Lemma
1 are true for $L$ and moreover we may replace everywhere $\phi ^2$ with $%
\phi ^3.$ With this in mind, let us fix $b_0=\delta \varepsilon ,$ $%
b_1=\beta T,$ where $0<\delta \varepsilon <\varepsilon <T<\beta T$ and for
convenience choose $\beta =\frac 1\delta $. Hence, everywhere in this proof 
\[
L=L\left( \delta \varepsilon ,\beta T\right) ,
\]
and $a=\frac m{m+1}.$

The proof of this theorem is based on a more detailed study of inequalities
in Lemma 1 and Lemma 2. We want to prove that $\lambda _1\left( M\right) =m^2
$ forces all the inequalities  to become equalities on $L\left( \varepsilon
,T\right) .$ Since $\varepsilon ,T$ are arbitrary, it will follow that we
need to have equalities everywhere on $M.$

Choose $\phi =\varphi \psi ,$ where 
\[
\psi =\left\{ 
\begin{array}{c}
1 \\ 
R-r \\ 
0
\end{array}
\left. 
\begin{array}{c}
\text{on} \\ 
\text{on} \\ 
\text{on}
\end{array}
\left. 
\begin{array}{l}
E_p\left( R-1\right) \cup F \\ 
E_p\left( R\right) \backslash E_p\left( R-1\right)  \\ 
E\backslash E_p\left( R\right) 
\end{array}
\right. \right. \right. 
\]
and 
\[
\varphi =\left\{ 
\begin{array}{c}
(-\log \delta )^{-1}(\log f-\log \left( \delta \varepsilon \right) ) \\ 
0 \\ 
(\log \beta )^{-1}(\log (\beta T)-\log f) \\ 
1
\end{array}
\left. 
\begin{array}{l}
\text{on}\;\;\;L\left( \delta \varepsilon ,\varepsilon \right)  \\ 
\text{on}\;\;\;L\left( 0, \delta \varepsilon \right)
\cup \left(L\left(\beta T, \infty \right)\cap F\right)  \\ 
\text{on}\;\;\;L(T,\beta T)\cap F \\ 
\text{otherwise.}
\end{array}
\right. \right. 
\]
Recall that by (\ref{6}) and (\ref{7}) we have 
\begin{gather}
\frac 1{16}(\frac 1m-\frac{\left( 1-a\right) ^2}{a(2-a)})\int_Lf^{-3}\left|
\nabla f\right| ^4\phi ^3\leq \frac a{2-a}\frac{m+1}4\int_Lf^{-1}\left|
\nabla f\right| ^2\phi ^3  \nonumber \\
+\frac 14\int_Lf^{-2}\left| \nabla f\right| ^2Re(f_{\bar{\beta}}\left( \phi
^3\right) _\beta )-\int_Lf^{-1}f_{\overline{\alpha }\beta }f_\alpha \left(
\phi ^3\right) _{\overline{\beta }}  \nonumber \\
-\frac a{2-a}\int_Lf^{-1}Re(f_\alpha f_{\bar{\alpha}\bar{\beta}}\left( \phi
^3\right) _\beta ).  \label{8}
\end{gather}
On the other hand, Schwarz inequality implies 
\begin{equation}
\left( \int_Lf^{-1}\left| \nabla f\right| ^2\phi ^3\right) ^2\leq \left(
\int_Lf^{-3}\left| \nabla f\right| ^4\phi ^3\right) \left( \int_Lf\phi
^3\right) ,  \label{9}
\end{equation}
and by the variational principle it follows 
\begin{eqnarray*}
\lambda _1\left( M\right) \int_Lf\phi ^3 &\leq &\int_L\left| \nabla \left(
f^{\frac 12}\phi ^{\frac 32}\right) \right| ^2 \\
&=&\frac 14\int_Lf^{-1}\left| \nabla f\right| ^2\phi ^3+\frac 94\int_Lf\phi
\left| \nabla \phi \right| ^2+\frac 32\int_L\phi ^2\nabla f\cdot \nabla \phi
.
\end{eqnarray*}
Our point now is that a careful study of the two $\nabla \phi -$terms shows
that they converge to zero as $\beta \rightarrow \infty $ (and $\delta
=\frac 1\beta \rightarrow 0$).

It is clear that $\frac 94\int_Lf\phi \left| \nabla \phi \right| ^2\leq 
\frac{c_1}{\log \beta },$ while

\[
\int_L\phi ^2\nabla f\cdot \nabla \phi =\frac 1{(-\log \delta
)}\int_{L\left( \delta \varepsilon ,\varepsilon \right) }f^{-1}\left| \nabla
f\right| ^2\phi ^2-\frac 1{\log \beta }\int_{L\left( T,\beta T\right) \cap
F}f^{-1}\left| \nabla f\right| ^2\phi ^2.
\]
The integral on $F$ is readily found by the co-area formula: 
\begin{eqnarray*}
\frac 1{\log \beta }\int_{L\left( T,\beta T\right) \cap F}f^{-1}\left|
\nabla f\right| ^2\phi ^2 &=&(\int_{l\left( t_0\right) }\left| \nabla
f\right| )\int_T^{\beta T}t^{-1}\left( \frac{\log (\beta T)-\log t}{\log
\beta }\right) ^2dt \\
&=&\frac 13\int_{l\left( t_0\right) }\left| \nabla f\right| .
\end{eqnarray*}
It is clear that the same formula holds on $E$ if we integrate against $%
\varphi ^2$ and therefore: 
\[
\frac 1{(-\log \delta )}\int_{L\left( \delta \varepsilon ,\varepsilon
\right) }f^{-1}\left| \nabla f\right| ^2\phi ^2\leq \frac 1{(-\log \delta
)}\int_{L\left( \delta \varepsilon ,\varepsilon \right) }\varphi ^2\left|
\nabla f\right| ^2f^{-1}=\frac 13\int_{l\left( t_0\right) }\left| \nabla
f\right| .
\]
For later use, observe that a converse of the latter inequality also holds: 
\begin{eqnarray*}
\frac 1{(-\log \delta )}\int_{L\left( \delta \varepsilon ,\varepsilon
\right) }f^{-1}\left| \nabla f\right| ^2\phi ^2 &\geq &\frac 1{(-\log \delta
)}\int_{L\left( \delta \varepsilon ,\varepsilon \right) }f^{-1}\left| \nabla
f\right| ^2\varphi ^2 \\
&&-\frac 1{(-\log \delta )}\int_{L\left( \delta \varepsilon ,\varepsilon
\right) \cap \left( E\backslash E_p\left( R-1\right) \right) }f^{-1}\left|
\nabla f\right| ^2\varphi ^2 \\
&\geq &\frac 13\int_{l\left( t_0\right) }\left| \nabla f\right| -\frac{c_2}{%
(-\log \delta )}.
\end{eqnarray*}
In particular, from the above estimates it follows 
\[
\int_L\phi ^2\nabla f\cdot \nabla \phi \leq 0.
\]
We have thus proved that 
\[
\lambda _1\left( M\right) \int_Lf\phi ^3\leq \frac 14\int_Lf^{-1}\left|
\nabla f\right| ^2\phi ^3+\frac{c_1}{\log \beta },
\]
which plugged into (\ref{9}) yields 
\begin{eqnarray*}
\int_Lf^{-3}\left| \nabla f\right| ^4\phi ^3 &\geq &4\lambda _1\left(
M\right) \frac{\left( \int_Lf^{-1}\left| \nabla f\right| ^2\phi ^3\right) ^2%
}{\int_Lf^{-1}\left| \nabla f\right| ^2\phi ^3+\frac{c_3}{\log \beta }} \\
&=&4\lambda _1\left( M\right) \int_Lf^{-1}\left| \nabla f\right| ^2\phi ^3-%
\frac{c_4}{\log \beta }\frac{\int_Lf^{-1}\left| \nabla f\right| ^2\phi ^3}{%
\int_Lf^{-1}\left| \nabla f\right| ^2\phi ^3+\frac{c_3}{\log \beta }} \\
&\geq &4\lambda _1\left( M\right) \int_Lf^{-1}\left| \nabla f\right| ^2\phi
^3-\frac{c_4}{\log \beta }.
\end{eqnarray*}
Now let's return to (\ref{8}) and use this lower bound, it follows that we
have 
\begin{eqnarray}
0 &\leq &\frac{c_5}{\log \beta }+\frac 14\int_Lf^{-2}\left| \nabla f\right|
^2Re(f_{\bar{\beta}}\left( \phi ^3\right) _\beta )  \nonumber \\
&&-\int_Lf^{-1}f_{\overline{\alpha }\beta }f_\alpha \left( \phi ^3\right) _{%
\overline{\beta }}-\frac a{2-a}\int_Lf^{-1}Re(f_\alpha f_{\bar{\alpha}\bar{%
\beta}}\left( \phi ^3\right) _\beta ).  \label{10}
\end{eqnarray}

\textbf{Claim:}

There exists a constant $c\geq 0$ such that

\begin{gather*}
\frac 14\int_Lf^{-2}\left| \nabla f\right| ^2Re(f_{\bar{\beta}}\left( \phi
^3\right) _\beta )-\int_Lf^{-1}f_{\overline{\alpha }\beta }f_\alpha \left(
\phi ^3\right) _{\overline{\beta }} \\
-\frac a{2-a}\int_Lf^{-1}Re(f_\alpha f_{\bar{\alpha}\bar{\beta}}\left( \phi
^3\right) _\beta )\leq \frac c{(\log \beta )^{\frac 12}}.
\end{gather*}
\textbf{Proof of the claim}.

Let us study each of the three terms in the left hand side.

I. We have: 
\begin{eqnarray*}
\frac 14\int_Lf^{-2}\left| \nabla f\right| ^2Re(f_{\bar{\beta}}\left( \phi
^3\right) _\beta ) &=&\frac 3{16}\int_L\phi ^2f^{-2}\left| \nabla f\right|
^2\nabla f\cdot \nabla \phi \\
&=&\frac 3{16}\frac 1{(-\log \delta )}\int_{L\left( \delta \varepsilon
,\varepsilon \right) }f^{-3}\left| \nabla f\right| ^4\phi ^2 \\
&&-\frac 3{16}\frac 1{\log \beta }\int_{L\left( T,\beta T\right) \cap
F}f^{-3}\left| \nabla f\right| ^4\phi ^2.
\end{eqnarray*}

As we stressed above, the estimates in Lemma 1 and Lemma 2 are true on any
end. Certainly, $\phi $ here is not the same on $L\left( \delta \varepsilon
,\varepsilon \right) $ with $\phi $ from Lemma 1. Nevertheless, the
computations are the same. In fact, in this case there is no need to
consider a cut-off $\varphi $ on $L\left( \frac 12\delta \varepsilon ,\delta
\varepsilon \right) ,$ because $\phi $ already is zero there. On $L\left(
\varepsilon ,2\varepsilon \right) $ the cut-off $\varphi $ is the same as in
Lemma 1. Therefore, one can use Lemma 1 for $L\left( \delta \varepsilon
,\varepsilon \right) $ and Lemma 2 for $L(T,\beta T)\cap F$ to estimate the
above subtraction.

By Lemma 1 we know that 
\begin{eqnarray*}
\frac 1{(-\log \delta )}\int_{L\left( \delta \varepsilon ,\varepsilon
\right) }f^{-3}\left| \nabla f\right| ^4\phi ^2 &\leq &4\lambda _1\left(
M\right) \frac 1{(-\log \delta )}\int_{L\left( \delta \varepsilon
,\varepsilon \right) }f^{-1}\left| \nabla f\right| ^2\phi ^2 \\
&&+\frac{c_6}{(-\log \delta )^{\frac 12}} \\
&\leq &\frac 43\lambda _1\left( M\right) \int_{l\left( t_0\right) }\left|
\nabla f\right| +\frac{c_6}{(-\log \delta )^{\frac 12}},
\end{eqnarray*}
while Lemma 2 implies that 
\[
\frac 1{\log \beta }\int_{L\left( T,\beta T\right) \cap F}f^{-3}\left|
\nabla f\right| ^4\phi ^2\geq \frac 43\lambda _1\left( M\right)
\int_{l\left( t_0\right) }\left| \nabla f\right| -\frac{c_7}{(\log \beta
)^{\frac 12}}. 
\]
Combining the two estimates, it results 
\[
\frac 14\int_Lf^{-2}\left| \nabla f\right| ^2Re(f_{\bar{\beta}}\left( \phi
^3\right) _\beta )\leq \frac{c_8}{(\log \beta )^{\frac 12}}. 
\]

II. Start with 
\begin{eqnarray*}
-\int_Lf^{-1}f_{\overline{\alpha }\beta }f_\alpha \left( \phi ^3\right) _{%
\overline{\beta }} &=&-\frac 3{(-\log \delta )}\int_{L\left( \delta
\varepsilon ,\varepsilon \right) }f^{-2}(f_{\overline{\alpha }\beta
}f_\alpha f_{\bar{\beta}})\phi ^2 \\
&&+\frac 3{\log \beta }\int_{L\left( T,\beta T\right) \cap F}f^{-2}(f_{%
\overline{\alpha }\beta }f_\alpha f_{\bar{\beta}})\phi ^2.
\end{eqnarray*}
From (\ref{1}) and (\ref{7}) we have: 
\begin{gather*}
\frac 1{\log \beta }\int_{L\left( T,\beta T\right) \cap F}f^{-2}(f_{%
\overline{\alpha }\beta }f_\alpha f_{\bar{\beta}})\phi ^2\leq \left( \frac
1{16}-\frac 1{16m}\right) \frac 1{\log \beta }\times  \\
\times \int_{L\left( T,\beta T\right) \cap F}f^{-3}\left| \nabla f\right|
^4\phi ^2+\frac{c_9}{(\log \beta )^{\frac 12}} \\
\leq\left( \frac 1{16}-\frac 1{16m}\right) \frac 43\lambda _1\left( M\right)
\int_{l\left( t_0\right) }\left| \nabla f\right| +\frac{c_{10}}{(\log \beta
)^{\frac 12}},
\end{gather*}
while from (\ref{5}) we know that 
\begin{gather*}
-\frac 1{(-\log \delta )}\int_{L\left( \delta \varepsilon ,\varepsilon
\right) }f^{-2}(f_{\overline{\alpha }\beta }f_\alpha f_{\bar{\beta}})\phi
^2\leq \left( -\frac 18+\frac 1{16a\left( 2-a\right) }\right) \frac 1{(-\log
\delta )}\times  \\
\times \int_{L\left( \delta \varepsilon ,\varepsilon \right) }f^{-3}\left|
\nabla f\right| ^4\phi ^2 \\
+\frac a{2-a}\frac{m+1}4\frac 1{(-\log \delta )}\int_{L\left( \delta
\varepsilon ,\varepsilon \right) }f^{-1}\left| \nabla f\right| ^2\phi ^2+%
\frac{c_{11}}{(-\log \delta )^{\frac 12}} \\
\leq \left( \left( -\frac 18+\frac 1{16a\left( 2-a\right) }\right) \frac
43\lambda _1\left( M\right) +\frac a{2-a}\frac{m+1}4\frac 13\right)
\int_{l\left( t_0\right) }\left| \nabla f\right|  \\
+\frac{c_{12}}{(-\log \delta )^{\frac 12}},
\end{gather*}
using the estimates in I. It is easy to see that the coefficients of $%
\int_{l\left( t_0\right) }\left| \nabla f\right| $ cancel out (this comes as
no surprise) and therefore 
\[
-\int_Lf^{-1}f_{\overline{\alpha }\beta }f_\alpha \left( \phi ^3\right) _{%
\overline{\beta }}\leq \frac{c_{13}}{(\log \beta )^{\frac 12}}.
\]
Note also that in a similar fashion it can be proved that 
\[
\int_Lf^{-1}f_{\overline{\alpha }\beta }f_\alpha \left( \phi ^3\right) _{%
\overline{\beta }}\leq \frac{c_{14}}{(\log \beta )^{\frac 12}}.
\]

III. Finally, by (\ref{2}) one has:

\begin{gather*}
-\int_Lf^{-1}Re(f_\alpha f_{\bar{\alpha}\bar{\beta}}\left( \phi ^3\right)
_\beta )=-\frac 3{(-\log \delta )}\int_{L\left( \delta \varepsilon
,\varepsilon \right) }f^{-2}Re(f_{\overline{\alpha }\bar{\beta}}f_\alpha f_{%
\bar{\beta}})\phi ^2 \\
+\frac 3{\log \beta }\int_{L\left( T,\beta T\right) \cap F}f^{-2}Re(f_{%
\overline{\alpha }\beta }f_\alpha f_{\bar{\beta}})\phi ^2 \\
\leq -\frac 3{8(-\log \delta )}\int_{L\left( \delta \varepsilon ,\varepsilon
\right) }f^{-3}\left| \nabla f\right| ^4\phi ^2+\frac 3{(-\log \delta
)}\int_{L\left( \delta \varepsilon ,\varepsilon \right) }f^{-2}(f_{\alpha 
\bar{\beta}}f_{\bar{\alpha}}f_\beta )\phi ^2 \\
+\frac 3{8\log \beta }\int_{L\left( T,\beta T\right) \cap F}f^{-3}\left|
\nabla f\right| ^4\phi ^2-\frac 3{\log \beta }\int_{L\left( T,\beta T\right)
\cap F}f^{-2}(f_{\alpha \bar{\beta}}f_{\bar{\alpha}}f_\beta )\phi ^2 \\
+\frac{c_{15}}{(\log \beta )^{\frac 12}}.
\end{gather*}
By I and II it can be proved that 
\[
\int_Lf^{-1}Re(f_\alpha f_{\bar{\alpha}\bar{\beta}}\left( \phi ^3\right)
_\beta )\leq \frac{c_{16}}{(\log \beta )^{\frac 12}}.
\]
This proves the claim. \textbf{Q.E.D.}

Let us use this result in (\ref{10}), then we infer that 
\begin{gather}
0\leq \frac{c_5}{\log \beta }+\frac 14\int_Lf^{-2}\left| \nabla f\right|
^2Re(f_{\bar{\beta}}\left( \phi ^3\right) _\beta )  \nonumber \\
-\int_Lf^{-1}f_{\overline{\alpha }\beta }f_\alpha \left( \phi ^3\right) _{%
\overline{\beta }}-\frac a{2-a}\int_Lf^{-1}Re(f_\alpha f_{\bar{\alpha}\bar{%
\beta}}\left( \phi ^3\right) _\beta )\leq \frac C{(\log \beta )^{\frac 12}}.
\label{11}
\end{gather}
Since $\beta $ (and $\delta =\frac 1\beta $) is arbitrary it follows that
for $\varepsilon $ and $T$ fixed the above inequality becomes equality by
letting $\beta \rightarrow \infty .$

From (\ref{11}) we are able to draw the conclusion that the following
formulas need to hold on $M$: 
\begin{eqnarray}
Ric_{1\bar{1}} &=&-(m+1)  \nonumber \\
\left| \nabla f\right|  &=&2\sqrt{\lambda _1\left( M\right) }f  \label{12} \\
u_{\alpha \bar{\beta}} &=&-m\delta _{\alpha \bar{\beta}}  \nonumber \\
u_{\alpha \beta } &=&m\delta _{1\alpha }\delta _{1\beta }  \nonumber
\end{eqnarray}
with respect to the frame 
\begin{eqnarray*}
v_\alpha  &=&\frac 12\left( e_{2\alpha -1}-\sqrt{-1}Je_{2\alpha -1}\right) ,
\\
e_1 &=&\frac 1{\left| \nabla f\right| }\nabla f,\;\;Je_{2k-1}=e_{2k}.
\end{eqnarray*}
Note that in view of (\ref{12}) this frame is globally defined on $M.$

Let us prove that indeed we have these relations on $M.$

Suppose that there exists a point $x_0\in M$ and a positive $\eta _0$ such
that: 
\[
Ric_{1\bar{1}}(x_0)\geq -(m+1)+\eta _0. 
\]
Let us choose $\varepsilon $ and $T$ such that $x_0\in L\left( \varepsilon
,T\right) .$ Recall that $L=L\left( \delta \varepsilon ,\beta T\right) ,$
for arbitrary $\beta $ and for $\delta =\frac 1\beta .$

Then one can see that there exists $\eta _1>0$ such that 
\[
-\int_Lf^{-1}f_\alpha f_{\bar{\alpha}\bar{\beta}\beta }\phi ^3\leq \frac{m+1}%
4\int_Lf^{-1}\left| \nabla f\right| ^2\phi ^3-\eta _1.
\]
It is easy to check that now (\ref{11}) will become

\begin{eqnarray*}
0 &<&\eta _1\leq \frac{c_5}{\log \beta }+\frac 14\int_Lf^{-2}\left| \nabla
f\right| ^2Re(f_{\bar{\beta}}\left( \phi ^3\right) _\beta ) \\
&&-\int_Lf^{-1}f_{\overline{\alpha }\beta }f_\alpha \left( \phi ^3\right) _{%
\overline{\beta }}-\frac a{2-a}\int_Lf^{-1}Re(f_\alpha f_{\bar{\alpha}\bar{%
\beta}}\left( \phi ^3\right) _\beta )\leq \frac C{(\log \beta )^{\frac 12}},
\end{eqnarray*}
which gives a contradiction if we let $\beta \rightarrow \infty .$

Next, let us focus on the Schwarz inequality (\ref{9}). Suppose by absurd
that there exists no constant $a\neq 0$ such that 
\[
\left| \nabla f\right| \left( x\right) =af\left( x\right) \;\text{for any}%
\;x\in U,
\]
where $U\subset L\left( \varepsilon ,T\right) $ is a fixed open set. It is
clear that if 
\[
h=f^{-\frac 32}\left| \nabla f\right| ^2\phi ^{\frac 32},\;\;g=f^{\frac
12}\phi ^{\frac 32},
\]
then there exists no $a\in \mathbb{R}$ such that $g=ah$ on $U,$ which implies
that 
\[
\eta _0:=\min_{a\in \mathbb{R}}\int_U\left( g-ah\right) ^2>0.
\]
This shows that 
\begin{eqnarray*}
\eta _0 &\leq &a^2\int_Uh^2-2a\int_Ugh+\int_Ug^2, \\
0 &\leq &a^2\int_{L\backslash U}h^2-2a\int_{L\backslash
U}gh+\int_{L\backslash U}g^2,
\end{eqnarray*}
for any $a\in \mathbb{R}.$ As a consequence, the following inequality is true
for any $a\in \mathbb{R}:$ 
\[
0\leq a^2\int_Lh^2-2a\int_Lgh+\left( \int_Lg^2-\eta _0\right) .
\]
It follows that 
\[
\left( \int_Lgh\right) ^2\leq \left( \int_Lh^2\right) \left( \int_Lg^2-\eta
_0\right) .
\]
Similarly, one can see that there exists an $\eta _1>0$ such that 
\[
\left( \int_Lgh\right) ^2\leq \left( \int_Lg^2\right) \left( \int_Lh^2-\eta
_1\right) .
\]
Adding these two inequalities and using the arithmetic mean inequality we
get that there exists $\eta _2>0$ with the property 
\[
\left( \int_Lgh+\eta _2\right) ^2\leq \int_Lg^2\int_Lh^2.
\]
We have thus proved that there exists a constant $\eta _2>0$ depending on $U$
but not on $\beta $ (and $\delta $) such that 
\[
\left( \int_Lf^{-1}\left| \nabla f\right| ^2\phi ^3+\eta _2\right) ^2\leq
\left( \int_Lf^{-3}\left| \nabla f\right| ^4\phi ^3\right) \left(
\int_Lf\phi ^3\right) ,
\]
inequality that will be used instead of (\ref{9}) in the argument that
followed. Consequently, 
\begin{eqnarray*}
\int_Lf^{-3}\left| \nabla f\right| ^4\phi ^3 &\geq &4\lambda _1\left(
M\right) \frac{\left( \int_Lf^{-1}\left| \nabla f\right| ^2\phi ^3+\eta
_2\right) ^2}{\int_Lf^{-1}\left| \nabla f\right| ^2\phi ^3+\frac{c_3}{\log
\beta }} \\
&\geq &4\lambda _1\left( M\right) \int_Lf^{-1}\left| \nabla f\right| ^2\phi
^3+8\lambda _1\left( M\right) \eta _2-c_{17}\frac 1{\log \beta }.
\end{eqnarray*}
However, using the same reasoning as for Ricci one can see that this yields
a contradiction.

Summing up, we have proved that there exists a constant $a>0$ such that $%
\left| \nabla f\right| =af$ on $M.$ Using Lemma 1 and Lemma 2 one can see
that $a=2\sqrt{\lambda _1\left( M\right) }.$

The proofs for the remaining two formulas use the same ideas. Note that in (%
\ref{7}) we need to have equality everywhere on $M,$ therefore there exists
a function $\mu $ on \thinspace $M$ such that 
\[
u_{\alpha \bar{\beta}}=\mu \delta _{\alpha \bar{\beta}}.
\]
However, taking the trace and using that $f$ is harmonic one can show that $%
\mu =-m.$

Finally, we pointed out that if equality holds in (\ref{3}) then 
\[
f_{\alpha \beta }=0\;\;\text{for}\;(\alpha ,\beta )\neq \left( 1,1\right) ,
\]
and on the other hand equality holds in (\ref{4}) if and only if 
\begin{eqnarray*}
\left| f_{\alpha \beta }\right|  &=&\frac 1a\frac{\left| \nabla f\right| ^2}{%
4f}=m(m+1)f, \\
Re(f_{11}) &=&\left| f_{11}\right| .
\end{eqnarray*}
This means that 
\[
f_{11}=m(m+1)f,
\]
or in terms of $u$ one has 
\[
u_{\alpha \beta }=m\delta _{1\alpha }\delta _{1\beta }
\]
as claimed.

Now we are ready to complete the proof of Theorem 4. 
Let us compute the real
Hessian of 
\[
B:=\frac 1{2m}u.
\]
We have: 
\begin{eqnarray*}
B_{e_1e_1} &=&B_{11}+B_{\bar{1}\bar{1}}+2B_{1\bar{1}}=1-1=0, \\
B_{e_2e_2} &=&-\left( B_{11}+B_{\bar{1}\bar{1}}-2B_{1\bar{1}}\right) =-2, \\
B_{e_{2k-1}e_{2k-1}} &=&B_{kk}+B_{\bar{k}\bar{k}}+2B_{k\bar{k}}=-1, \\
B_{e_{2k}e_{2k}} &=&-B_{kk}-B_{\bar{k}\bar{k}}+2B_{k\bar{k}}=-1, \\
B_{e_ke_j} &=&0\;\;if\;\;k\neq j,
\end{eqnarray*}
for $k\in \{2,...,m\}.$ Also, notice that $\left| \nabla B\right| =1$ on $M.$

Since all the computations from now on will be done in the real frame $%
\left\{ e_1,....,e_{2m}\right\} $ with $Je_{2k-1}=e_{2k}$ and $e_1=\frac
1{\left| \nabla f\right| }\nabla f,$ for convenience we will drop the $e_k$
index and use only $k$ in the formulas for the real Hessian and the
curvature.

Also, let us make the convention that Roman letters $i,j,k$ run from $1$ to $%
2m$ and Greek letters $\alpha ,\beta ,\gamma $ run from $3$ to $2m.$

We have proved that there exists a smooth function $B$ on $M$ with real
Hessian 
\[
\left( B_{ij}\right) =\left( 
\begin{array}{lllllll}
0 & 0 & 0 & 0 & . & . & 0 \\ 
0 & -2 & 0 & 0 & . & . & 0 \\ 
0 & 0 & -1 & 0 & . & . & 0 \\ 
0 & 0 & 0 & -1 & . & . & 0 \\ 
0 & . & . & . & . & . & . \\ 
0 & . & . & . & . & . & . \\ 
0 & 0 & 0 & 0 & . & . & -1
\end{array}
\right) 
\]
and with unit gradient, $\left| \nabla B\right| =1$ on $M.$

Note that our function $B$ satisfies the same properties as the Buseman function 
$\beta$ in \cite{L-W}. Using the Hessian of $\beta$ Li-Wang proved that the manifold 
has at most two ends and in the two ends case they inferred some information on the structure 
of $M$. 
We will give an outline of their argument below.

Denote the level set of $B$ by 
\[
N_t=\left\{ x\in M\left| \;B\left( x\right) =t\right. \right\} .
\]

Since $\left| \nabla B\right| =1$, $M$ is diffeomorphic to $\mathbb{R}\times N_0
$ and $e_1=\nabla B$ is the unit normal to $N_t$ for any $t.$ If $N_0$ is noncompact,
then $M$ will have one end, wich contradicts our assumption that $M$ has more than one end.

Consequently, $N_0$ is compact, and this implies that $M$ has two ends. 
For the remainder of this proof $M$ has two ends, and we want to find the metric of $N_t$ 
depending on the metric of $N_0$.
 
Knowing $B_{ij}$ is equivalent to knowing the second fundamental form of $%
N_t,$ which implies that if 
\[
\nabla e_i=\omega _{ik}e_k,
\]
then one can find 
\[
\omega _{i1}\left( e_j\right) =\left\{ 
\begin{array}{l}
0\;\;\;\text{for}\;\;i\neq j \\ 
2\;\;\;\text{for\ }\;i=j=2 \\ 
1\;\;\;\text{for}\;\;3\leq i=j\leq 2m.
\end{array}
\right. 
\]
Also, using the K\"{a}hler property we know that 
\[
\omega _{1k}Je_k=J\nabla e_1=\nabla Je_1=\nabla e_2=\omega _{2k}e_k,
\]
which implies 
\[
\omega _{\alpha 2}\left( e_j\right) =\left\{ 
\begin{array}{l}
\;0\;\;\;\;\text{for}\;\;j=1\;or\;j=2 \\ 
-1\;\;\;\text{for}\;\;\alpha =2p+1,\;j=2p+2\;\; \\ 
\;1\;\;\;\text{for\ }\;\alpha =2p+2,\;j=2p+1.
\end{array}
\right. 
\]
It is clear that the flow $\phi _t:M\rightarrow M$ generated by $e_1$ is a
geodesic flow. Since 
\[
\nabla _{e_1}e_2=\nabla _{e_1}Je_1=J\nabla _{e_1}e_1=0
\]
we can conclude that $e_2$ is parallel along the geodesic $\tau $ defined by 
$e_1$. We will consider the rest of the frame so that it is also parallel
along this geodesic.

The next step is to prove that 
\begin{eqnarray*}
V_2\left( t\right)  &=&e^{-2t}e_2 \\
V_\alpha \left( t\right)  &=&e^{-t}e_\alpha 
\end{eqnarray*}
are the Jacobi fields along the geodesic $\tau $ with initial conditions 
\begin{eqnarray*}
V_2\left( 0\right)  &=&e_2,\;\;V_2^{\prime }\left( 0\right) =-2e_2 \\
V_\alpha \left( 0\right)  &=&e_\alpha ,\;\;V_\alpha ^{\prime }\left(
0\right) =-e_\alpha .
\end{eqnarray*}
This is true because the information on $\omega _{i1}$ and $\omega _{\alpha
2}$ allows to find sufficient values for the curvature tensor. Using the
second structural equations one can show that 
\begin{eqnarray}
R_{1212} &=&-4,\;\;R_{121\alpha }=0,  \nonumber \\
R_{1\alpha 1\beta } &=&-\delta _{\alpha \beta },  \label{13}
\end{eqnarray}
and this indeed shows that $V_k\left( t\right) $ are Jacobi fields for $k\in
\{2,...,2m\}.$

However, $d\phi _t\left( e_k\right) $ for $k\geq 2$ are also Jacobi fields
with the same initial conditions as $V_k\left( t\right) ,$ so they must
coincide.

The conclusion is that the metrics on $N_t$ viewed as one parameter of
metrics on $N_0$ are 
\[
ds_t^2=e^{-4t}\omega _2^2\left( 0\right) +e^{-2t}\left( \omega _3^2\left(
0\right) +....+\omega _{2m}^2\left( 0\right) \right) , 
\]
where $\left\{ \omega _1,...,\omega _{2m}\right\} $ is the dual frame of $%
\left\{ e_1,...,e_{2m}\right\} $.\textbf{\ Q.E.D.}

{\scriptsize DEPARTMENT OF MATHEMATICS, UNIVERSITY OF CALIFORNIA, IRVINE,
CA, 92697-3875}

{\small E-mail address: omuntean@math.uci.edu}

\end{document}